\documentclass[10pt,oneside]{article}
\input amssymb.sty
\usepackage[cp1250]{inputenc}
\usepackage[T1]{fontenc}

\begin{document}
\newcommand{\qed}{\rule{1.5mm}{1.5mm}}
\newcommand{\proof}{\textit{Proof. }}
\newcommand{\ccon}{\rightarrowtail}
\newtheorem{theorem}{Theorem}[section]
\newtheorem{lemma}[theorem]{Lemma}
\newtheorem{remark}[theorem]{Remark}
\newtheorem{example}[theorem]{Example}
\newtheorem{corollary}[theorem]{Corollary}
\newtheorem{proposition}[theorem]{Proposition}
\newtheorem{claim}[theorem]{Claim}
\begin{center}
{\LARGE\textbf{Algebraic approximation of analytic sets and mappings}\vspace*{3mm}}\\
{\large Marcin Bilski\footnote{M. Bilski: Institute of
Mathematics, Jagiellonian University, Reymonta 4, 30-059 Krak\'ow,
Poland. e-mail: Marcin.Bilski@im.uj.edu.pl\vspace*{1mm}\\
This is a revised version of the paper submitted for publication
to a journal (receipt acknow\-ledged on 30 May 2007). Manuscript
revised on: 20 June 2007, 01 December 2007.\vspace*{1mm}\\
Research partially supported by the Polish Ministry of Science and
Higher Education.
}}\vspace*{8mm}\\
\end{center}
\begin{abstract}Let $\{X_{\nu}\}$ be a sequence of analytic sets
converging to some analytic set $X$ in the sense of holomorphic
chains. We introduce a condition which implies that every
irreducible component of $X$ is the limit of a sequence of
irreducible components of the sets from $\{X_{\nu}\}.$ Next we
apply the condition to approximate a holomorphic solution $y=f(x)$
of a system $Q(x,y)=0$ of Nash equations by Nash solutions.
Presented methods allow to construct an algorithm of
approximation of the holomorphic solutions.\vspace*{2mm}\vspace*{0mm}
\end{abstract}

\section{Introduction}

Let $\mathbf{K}$ denote the field of complex or real numbers. The
following approximation theorem is known to be true: every
$\mathbf{K}$-analytic mapping $f:\Omega\rightarrow\mathbf{K}^k$
such that $Q(x,f(x))=0$ for $x\in\Omega,$ where $Q$ is a
$\mathbf{K}$-Nash mapping ($\Omega$ described below), can be
uniformly approximated by a $\mathbf{K}$-Nash mapping
$F:\Omega\rightarrow\mathbf{K}^k$ such that $Q(x,F(x))=0$ for
$x\in\Omega.$

In the complex case the theorem was proved by L. Lempert (see
\cite{Lem}, Theorem~3.2) for every holomorphically convex compact
subset $\Omega$ of an affine\linebreak algebraic variety and in
the real case it was proved by M. Coste, J. Ruiz\linebreak and
M.~Shiota (see \cite{CRS}, Theorem 1.1) for every compact Nash
manifold $\Omega$ . The\linebreak approximation theorem turned out
to be a very strong tool with many important applications (see
\cite{CRS}, \cite{Lem}).

The proofs of the theorem presented in \cite{CRS}, \cite{Lem} rely
on the solution to the M.~Artin's conjecture: a deep result from
commutative algebra for which the reader is referred to \cite{An},
\cite{Og}, \cite{Po1}, \cite{Po2}, \cite{Sp}. Such an approach
enabled to reach the goal in an elegant and relatively short way.
On the other hand, it seems to be very difficult to apply the
proofs in order to find Nash approximations for concrete analytic
mappings hence it is natural to ask whether the theorem can be
obtained directly. The latter question is strongly motivated by
the fact that approximating analytic objects by algebraic
counterparts is one of central techniques used in numerical
computations. From this point of view it is important to develop
theory of approximation that could serve as a base for finding
numerical algorithms.

In Section \ref{mainproof} of the present paper we give, using
only some basic methods of analytic geometry, a proof of a
semi-global version of the theorem in the complex case (see
Theorem \ref{semiglobmain}). The proof allows to construct an
algorithm of approximation of the mapping $f$ which is described
in Subsection \ref{examp}.

The following local version is an immediate consequence of Theorem
\ref{semiglobmain}.
\begin{theorem}\label{main}Let $U$ be an open subset of
$\mathbf{C}^n$ and let $f:U\rightarrow\mathbf{C}^k$ be a
holomorphic mapping that satisfies a system of equations
$Q(x,f(x))=0$ for $x\in U.$ Here $Q$ is a Nash mapping from a
neighborhood $\hat{U}$ in $\mathbf{C}^n\times\mathbf{C}^k$ of the
graph of $f$ into some $\mathbf{C}^q.$ Then for every $x_0\in U $
there are an open neighborhood $U_0\subset U$ and a sequence
$\{f^{\nu}:U_0\rightarrow\mathbf{C}^k\}$ of Nash mappings
converging uniformly to $f|_{U_0}$ such that $Q(x,f^{\nu}(x))=0$
for every $x\in U_0$ and $\nu\in\mathbf{N}.$
\end{theorem}
In the local situation the problem of approximation of the
solutions of algebraic or analytic equations was investigated by
M.~Artin in \cite{Ar1}, \cite{Ar2}, \cite{Ar3} and
Theorem~\ref{main} can be derived from his results.

Our interest in Theorem \ref{main} and its generalizations is
partially motivated by applications in the theory of analytic
sets. In particular, papers \cite{B1}--\cite{B4} contain a variety
of results on approximation of complex analytic sets by complex
Nash sets whose proofs can be divided into two stages: (i)
preparation, where only direct geometric methods appear, (ii)
switching Theorem \ref{main}. Thus the techniques of the present
article allow to obtain many of these results in a purely
geometric way. As an example let us mention the following main
theorem of \cite{B5}. Let $X$ be an analytic subset of pure
dimension $n$ of an open set $U\subset\mathbf{C}^m$ and let $E$ be
a Nash subset of $U$ such that $E\subset X.$ Then for every $a\in
E$ there is an open neighborhood $U_a$ of $a$ in $U$ and a
sequence $\{X_{\nu}\}$ of complex Nash subsets of $U_a$ of pure
dimension $n$ converging to $X\cap U_a$ in the sense of
holomorphic chains such that the following hold for every
$\nu\in\mathbf{N}:$ $E\cap U_a\subset X_{\nu}$ and
$\mu_x(X_{\nu})=\mu_x(X)$ for $x\in (E\cap U_a)\setminus F_{\nu},$
where $F_{\nu}$ is a nowhere dense analytic subset of $E\cap U_a.$
Here $\mu_x(X)$ denotes the multiplicity of $X$ at $x$ (see
\cite{Ch}, \cite{D} for the properties and generalizations of this
notion).

In the proof of Theorem \ref{semiglobmain} we apply
Theorem~\ref{component} from Section \ref{firstmain} which, being
of independent interest, is the first main result of this paper.
The aim of Section \ref{firstmain} is to develop a method of
controlling the behavior of irreducible components of analytic
sets from a sequence $\{X_{\nu}\}$ converging in the sense of
holomorphic chains to some analytic set $X.$ More precisely, we
formulate conditions which guarantee that the numbers of the
irreducible components of $X$ and of $X_{\nu}$ are equal for
almost all $\nu$ which in the considered context implies that
every irreducible component of $X$ is the limit of a sequence of
irreducible components of the sets from $\{X_{\nu}\}.$

Combining (the global version of) Theorem \ref{main} with
Theorem~\ref{component} one obtains a new method of algebraic
approximation of analytic sets extending the approach of
\cite{B1}. Namely, let $X$ be an analytic subset of
$U\times\mathbf{C}^k$ of pure dimension $n$ with proper projection
onto the Runge domain $U\subset\mathbf{C}^n.$ It is well known
(\cite{Wh}) that $X$ is a subset of another purely $n$-dimensional
analytic set $X'$ given by
$$X'=\{(x,z_1,\ldots,z_k)\in U\times\mathbf{C}^k:
p_1(x,z_1)=\ldots=p_k(x,z_k)=0\},$$ where
$p_i\in\mathcal{O}(U)[z_i]$ is unitary with non-zero discriminant,
for $i=1,\ldots,k.$ After replacing the coefficients of $p_i,$ for
every $i,$ by their Nash approxi\-mations on $U$ we obtain the set
$\tilde{X}'$ approximating $X'.$ Clearly, this does not mean that
some components of $\tilde{X}'$ automatically approximate $X.$
Yet, by Theorem~\ref{component} there is a system of polynomial
equations satisfied by the coefficients of $p_i,$ $i=1,\ldots,k,$
with the following property. Let $\tilde{U}$ be any open
relatively compact subset of $U.$ If the Nash approximations of
the coefficients (used to define $\tilde{X}'$) are close enough to
the original coefficients and also satisfy the equations then
$X\cap(\tilde{U}\times\mathbf{C}^k)$ is approximated by some
components of $\tilde{X}'\cap(\tilde{U}\times\mathbf{C}^k).$ The
existence of the Nash approximations of the coefficients
satisfying the equations mentioned above in a neighborhood of
$\overline{\tilde{U}}$ follows by the global version of Theorem
\ref{main}.

Finally let us recall that the convergence of a sequence of
analytic sets in the sense of chains, appearing in Theorem
\ref{component}, is equivalent to the (introduced in \cite{Lel})
convergence of currents of integration over these sets in the
weak-$\star$ topology. The basic facts on holomorphic chains (and
preliminaries on Nash sets and analytic sets with proper
projection) are gathered in Section \ref{preli} below.
\section{Preliminaries}\label{preli}
\subsection{Nash sets}\label{prelnash}
Let $\Omega$ be an open subset of $\mathbf{C}^n$ and let $f$ be a
holomorphic function on $\Omega.$ We say that $f$ is a Nash
function at $x_0\in\Omega$ if there exist an open neighborhood $U$
of $x_0$ and a polynomial
$P:\mathbf{C}^n\times\mathbf{C}\rightarrow\mathbf{C},$ $P\neq 0,$
such that $P(x,f(x))=0$ for $x\in U.$ A holomorphic function
defined on $\Omega$ is said to be a Nash function if it is a Nash
function at every point of $\Omega.$ A holomorphic mapping defined
on $\Omega$ with values in $\mathbf{C}^N$ is said to be a Nash
mapping if each of its components is a Nash function.

A subset $Y$ of an open set $\Omega\subset\mathbf{C}^n$ is said to
be a Nash subset of $\Omega$ if and only if for every
$y_0\in\Omega$ there exists a neighborhood $U$ of $y_0$ in
$\Omega$ and there exist Nash functions $f_1,\ldots,f_s$ on $U$
such that $$Y\cap U=\{x\in U: f_1(x)=\ldots=f_s(x)=0\}.$$ We will
use the following fact from \cite{Tw}, p. 239. Let
$\pi:\Omega\times\mathbf{C}^k\rightarrow\Omega$ denote the natural
projection.
\begin{theorem}
\label{projnash} Let $X$ be a Nash subset of
$\Omega\times\mathbf{C}^k$ such that $\pi|_{X}:X\rightarrow\Omega$
is a proper mapping. Then $\pi(X)$ is a Nash subset of $\Omega$
and $dim(X)=dim(\pi(X)).$
\end{theorem}
The fact from \cite{Tw} stated below explains the relation between
Nash and algebraic sets.
\begin{theorem}Let $X$ be a Nash subset of an open set
$\Omega\subset\mathbf{C}^n.$ Then every analytic irreducible
component of $X$ is an irreducible Nash subset of $\Omega.$
Moreover, if $X$ is irreducible then there exists an algebraic
subset $Y$ of $\mathbf{C}^n$ such that $X$ is an analytic
irreducible component of $Y\cap\Omega.$
\end{theorem}

\subsection{Analytic sets}\label{setswithprop}
Let $U, U'$ be domains in $\mathbf{C}^n,\mathbf{C}^k$ respectively
and let
$\pi:\mathbf{C}^n\times\mathbf{C}^k\rightarrow\mathbf{C}^n$ denote
the natural projection. For any purely $n$-dimensional analytic
subset $Y$ of $U\times U'$ with proper projection onto $U$ by
$\mathcal{S}(Y,\pi)$ we denote the set of singular points of
$\pi|_{Y}:$
$$\mathcal{S}(Y,\pi)=Sing(Y)\cup\{x\in Reg(Y):(\pi|_Y)'(x) \mbox{ is not an isomorphism}\}.$$
We often write $\mathcal{S}(Y)$ instead of $\mathcal{S}(Y,\pi)$
when it is clear which projection is taken into consideration.

It is well known that $\mathcal{S}(Y)$ is an analytic subset of
$U\times U',$ $dim(Y)<n$ (cp. \cite{Ch}, p. 50), hence by the
Remmert theorem $\pi(\mathcal{S}(Y))$ is also analytic. Moreover,
the following hold. The mapping $\pi|_{Y}$ is surjective
and open and there exists an integer $s=s(\pi|_{Y})$ such that:\vspace*{2mm}\\
(1) $\sharp(\pi|_{Y})^{-1}(\{a\})<s$ for
$a\in\pi(\mathcal{S}(Y)),$\\
(2) $\sharp(\pi|_{Y})^{-1}(\{a\})=s$ for $a\in U\setminus
\pi(\mathcal{S}(Y)),$\\
(3) for every $a\in U\setminus \pi(\mathcal{S}(Y))$ there exists a
neighborhood $W$ of $a$ and holomor-\linebreak \hspace*{5.4mm}phic
mappings $f_1,\ldots,f_s:W\rightarrow U'$ such that $f_i\cap
f_j=\emptyset$ for $i\neq j$ and\linebreak
\hspace*{5.4mm}$f_1\cup\ldots\cup f_s=$$(W\times U')\cap
Y.$\vspace*{2mm}

Let $E$ be a purely $n$-dimensional analytic subset of $U\times
U'$ with proper projection onto a domain $U\subset\mathbf{C}^n,$
where $U'$ is a domain in $\mathbf{C}.$ The unitary polynomial
$p\in\mathcal{O}(U)[z]$ such that $E=\{(x,z)\in
U\times\mathbf{C}:p(x,z)=0\}$ and the discriminant of $p$ is not
identically zero will be called the optimal polynomial for $E.$

Finally, for any analytic subset $X$ of an open set
$\tilde{U}\subset\mathbf{C}^m$ let $X_{(k)}\subset\tilde{U}$
denote the union of all irreducible components of $X$ of dimension
$k.$
\subsection{Convergence of closed sets and holomorphic chains}\label{holchai}
Let $U$ be an open subset in $\mathbf{C}^m.$ By a holomorphic
chain in $U$ we mean the formal sum $A=\sum_{j\in J}\alpha_jC_j,$
where $\alpha_j\neq 0$ for $j\in J$ are integers and
$\{C_j\}_{j\in J}$ is a locally finite family of pairwise distinct
irreducible analytic subsets of $U$ (see \cite{Tw2}, cp. also
\cite{Ba}, \cite{Ch}). The set $\bigcup_{j\in J}C_j$ is called the
support of $A$ and is denoted by $|A|$ whereas the sets $C_j$ are
called the components of $A$ with multiplicities $\alpha_j.$ The
chain $A$ is called positive if $\alpha_j>0$ for all $j\in J.$ If
all the components of $A$ have the same dimension $n$ then $A$
will be called an $n-$chain.

Below we introduce the convergence of holomorphic chains in $U$.
To do this we first need the notion of the local uniform
convergence of closed sets (cp. \cite{TwW}). Let $Y,Y_{\nu}$ be
closed subsets of $U$ for $\nu\in\mathbf{N}.$ We say that
$\{Y_{\nu}\}$ converges to
$Y$ locally uniformly if:\vspace*{2mm}\\
(1l) for every $a\in Y$ there exists a sequence $\{a_{\nu}\}$ such
that $a_{\nu}\in Y_{\nu}$ and $a_{\nu}\rightarrow a$\linebreak
\hspace*{6.6mm}in the standard topology of $\mathbf{C}^m,$\\
(2l) for every compact subset $K$ of $U$ such that $K\cap
Y=\emptyset$ it holds $K\cap Y_{\nu}=\emptyset$\linebreak
\hspace*{6.6mm}for almost all $\nu.$\vspace*{1mm}\\
Then we write $Y_{\nu}\rightarrow Y.$

We say that a sequence $\{Z_{\nu}\}$ of positive $n-$chains
converges to a positive $n-$chain $Z$ if:\vspace*{2mm}\\
(1c) $|Z_{\nu}|\rightarrow |Z|,$\\
(2c) for each regular point $a$ of $|Z|$ and each submanifold $T$
of $U$ of dimension\linebreak \hspace*{7mm}$m-n$ transversal to
$|Z|$ at $a$ such that $\overline{T}$ is compact and
$|Z|\cap\overline{T}=\{a\},$\linebreak \hspace*{7mm}we have
$deg(Z_{\nu}\cdot T)=deg(Z\cdot T)$ for almost
all $\nu.$\vspace*{1mm}\\
Then we write $Z_{\nu}\ccon Z.$ By $Z\cdot T$ we denote the
intersection product of $Z$ and $T$ (cf. \cite{Tw2}). Observe that
the chains $Z_{\nu}\cdot T$ and $Z\cdot T$ for sufficiently large
$\nu$ have finite supports and the degrees are well defined.
Recall that for a chain $A=\sum_{j=1}^d\alpha_j\{a_j\},$
$deg(A)=\sum_{j=1}^d\alpha_j.$

The following lemma from \cite{Tw2} will be useful to us.
\begin{lemma}\label{eqconv} Let $n\in\mathbf{N}$ and $Z,Z_{\nu},$
for $\nu\in\mathbf{N},$ be positive $n$-chains. If
$|Z_{\nu}|\rightarrow |Z|$ then the following conditions are
equivalent:\\
(1) $Z_{\nu}\ccon Z,$\\
(2) for each point $a$ from a given dense subset of $Reg(|Z|)$
there is a submani-\linebreak\hspace*{5.9mm}fold $T$ of $U$ of
dimension $m-n$ transversal to $|Z|$ at $a$ such that
$\overline{T}$ is
compact,\linebreak\hspace*{5.7mm}$|Z|\cap\overline{T}=\{a\}$ and
$deg(Z_{\nu}\cdot T)=deg(Z\cdot T)$ for almost all $\nu.$
\end{lemma}

Let $U\subset\mathbf{C}^n$ be a domain and let
$\pi:U\times\mathbf{C}^k\rightarrow U$ be the natural projection.
Theorem \ref{manusc} below, taken from \cite{B1}, will be applied
in the proof of the main result ($s(\pi|_Y)$ is defined in Section
\ref{setswithprop}).
\begin{theorem}\label{manusc}Let $Y,Y_{\nu},$ for $\nu\in\mathbf{N},$
be purely $n$-dimensional analytic subsets of
$U\times\mathbf{C}^k$ with proper projection onto $U$ such that
$\{Y_{\nu}\}$ converges to $Y$ locally uniformly and
$s(\pi|_{Y_{\nu}})=s(\pi|_{Y})$ for every $\nu.$ Moreover, assume
that for every $\nu$ the number of the irreducible components of
$Y$ does not exceed the number of the irreducible components of
$Y_{\nu}.$ Then for each irreducible component $A$ of $Y$ there is
a sequence $\{A_{\nu}\}$ converging to $A$ locally uniformly such
that every $A_{\nu}$ is an irreducible component of $Y_{\nu}$ and
$s(\pi|_{A_{\nu}})=s(\pi|_A)$ for almost all $\nu.$
\end{theorem}

\section{Approximation}\label{proofs}
\subsection{Approximation of components of analytic sets}\label{firstmain}
Our first main result is the following theorem. Let
$U\subset\mathbf{C}^n$ be a domain and let
$\pi:U\times\mathbf{C}^k\rightarrow U$ denote the natural
projection. Let $X\subset U\times\mathbf{C}^k$ be an analytic
subset of pure dimension $n$ with proper projection onto $U.$
Recall that $s(\pi|_{X})$ denotes the cardinality of the generic
fiber in $X$ over $U.$
\begin{theorem}\label{component} Let $\{X_{\nu}\}$ be a sequence
of purely $n$-dimensional analytic subsets of
$U\times\mathbf{C}^k$ with proper projection onto $U$ converging
locally uniformly to $X$ such that $s(\pi|_{X})=s(\pi|_{X_{\nu}})$
for $\nu\in\mathbf{N}.$ Assume that
$\{(\pi(\mathcal{S}(X_{\nu})))_{(n-1)}\}$ converges to
$(\pi(\mathcal{S}(X)))_{(n-1)}$ in the sense of holomorphic
chains. Then for every analytic subset $Y$ of
$U\times\mathbf{C}^k$ of pure dimension $n$ such that $Y\subset X$
and for every open relatively compact subset $\tilde{U}$ of $U$
there exists a sequence $\{Y_{\nu}\}$ of purely $n$-dimensional
analytic subsets of $\tilde{U}\times\mathbf{C}^k$ converging to
$Y\cap(\tilde{U}\times\mathbf{C}^k)$ in the sense of holomorphic
chains such that $Y_{\nu}\subset X_{\nu}$ for every
$\nu\in\mathbf{N}.$
\end{theorem}
\begin{remark}\em From the proof of Theorem \ref{component} it follows
that if $Y\cap(\tilde{U}\times\mathbf{C}^k)$ is irreducible than
$Y_{\nu}$ is irreducible as well for almost all $\nu.$
\end{remark}
\textit{Proof of Theorem \ref{component}.} First the theorem will
be
proved under the following extra hypotheses:\vspace*{2mm}\\
(1) $U=U_1\times U_2\subset\mathbf{C}^{n-1}\times\mathbf{C}$ where
$U_1,$ $U_2$ are open balls,\\
(2) $\pi(\mathcal{S}(X))$ is with proper projection onto $U_1,$\\
(3) there is a compact ball $B$ in $U_2$ such that $(U_1\times
B)\cap \pi(\mathcal{S}(X))=\emptyset.$\vspace*{2mm}\\
This will be done in two steps. Step 1 is of preparatory nature.
Here we specify several technical conditions which may be assumed
satisfied by $X$ and $X_{\nu}$ (for large $\nu$) without loss of
generality. These conditions are used in Step 2 the idea of which
is to show that for almost all $\nu$ the number of the irreducible
components of $X\cap(\tilde{U}\times\mathbf{C}^k)$ does not exceed
the number of the irreducible components of
$X_{\nu}\cap(\tilde{U}\times\mathbf{C}^k)$ (in fact the numbers
are equal). This is done by constructing an injective mapping
which assigns to every irreducible component of
$X\cap(\tilde{U}\times\mathbf{C}^k)$ an irreducible component of
$X_{\nu}\cap(\tilde{U}\times\mathbf{C}^k).$ Then Theorem
\ref{manusc} from Section~\ref{holchai} may be applied which
completes the proof (under the hypotheses (1), (2), (3) above).

Finally we show (Step 3) that (1), (2), (3) are not necessary.\vspace*{2mm}\\
\textit{Step 1.} Without loss of generality we assume that $X$ and
$X_{\nu}$ (for large $\nu$) satisfy the conditions specified
below.

Denote
$$\hat{k}:=\max\{\sharp((\{x'\}\times U_2)\cap
(\pi(\mathcal{S}(X)))_{(n-1)}):x'\in U_1\}.$$ Let $\Sigma'(X)$ be
the subset of $U_1$ of points $x'$ for which
$$\sharp((\{x'\}\times U_2)\cap
(\pi(\mathcal{S}(X)))_{(n-1)})<\hat{k}.$$ Put
$\Sigma(X)=\Sigma'(X)\cup\rho(\overline{\pi(\mathcal{S}(X))\setminus\pi(\mathcal{S}(X))_{(n-1)}}),$
where $\rho:U_1\times U_2\rightarrow U_1$ is the natural
projection. (The closure is taken in $U_1\times U_2.$ Generally,
in this paper, the topological structure on any analytic set is
induced by the standard topology of $\mathbf{C}^m$ in which the
set is contained.)

Observe that $\Sigma(X)$ is a nowhere dense analytic subset of
$U_1$ hence there are $x'_0\in U_1\setminus\Sigma(X)$ and compact
balls $B_1,\ldots,B_{\hat{k}}\subset U_2$ such that
$B\cap(\bigcup_{i=1}^{\hat{k}}B_i)=\emptyset$ and $B_i\cap
B_j=\emptyset$ for $i\neq j.$ Moreover, each $B_i$ contains
precisely one $y$ such that
$(x'_0,y)\in\pi(\mathcal{S}(X))_{(n-1)}.$ Furthermore, since $U$
may be replaced by its relatively compact subset containing the
fixed $\tilde{U},$ we may assume that there is $r>0$ such that for
every $$x\in (U_1\times
B)\cup(\{x'_0\}\times(U_2\setminus(B_1\cup\ldots\cup
B_{\hat{k}})))$$ if $(x,u),(x,v)\in X$ and $u\neq v$ then
$||u-v||_{\mathbf{C}^k}>r.$

Next, by the fact that $\{\pi(\mathcal{S}(X_{\nu}))_{(n-1)}\}$
converges to $\pi(\mathcal{S}(X))_{(n-1)}$ in the sense of chains
and again by the fact that one may pass on to a relatively compact
subset of $U,$ the following is assumed for large $\nu:$
$\pi(\mathcal{S}(X_{\nu}))$ is with proper projection onto $U_1$
and the cardinality of the generic fiber of
$\pi(\mathcal{S}(X_{\nu}))_{(n-1)}$ over $U_1$ equals $\hat{k}.$
Moreover, in every $B_i$ there is precisely one $y$ such that
$(x'_0,y)\in\pi(\mathcal{S}(X_{\nu}))_{(n-1)}.$

Fix $x_0\in U_1\times B.$ Let $A\subset U_1\times
B\times\mathbf{C}^k$ denote the fiber in $X$ over $x_0.$ For every
irreducible component $Y$ of $X$ define $A_Y:=Y\cap A.$ For every
such $Y$ there is an arc
$$\gamma_Y:[0,1]\rightarrow ((U\setminus(\pi(\mathcal{S}(X))))\times\mathbf{C}^k)\cap Y$$
connecting all the points in $A_Y.$ Let
$$r_0=\inf\{||u-v||:(x,v),(x,u)\in X, u\neq v,
x\in\bigcup_{Y}\pi(\gamma_Y([0,1]))\}.$$ Then $r_0>0.$

Pick any $0<\delta<\min(\frac{r}{3},\frac{r_0}{3}).$ We complete
\textit{Step 1} by observing that for large $\nu$ the following
may be assumed:
$$(\bigcup_Y\pi(\gamma_{Y}([0,1])))\cap\pi(\mathcal{S}(X_{\nu}))=\emptyset$$ and
$$dist((\{x\}\times\mathbf{C}^k)\cap
X,(\{x\}\times\mathbf{C}^k)\cap X_{\nu})<\delta$$ for every $x\in
U$ (the latter due to the fact that $U$ may be replaced by its
relatively compact subset). (Here $dist$ denotes the Hausdorff
distance.)\vspace*{2mm}\\
\textit{Step 2.} We show that if $X$ and $X_{\nu}$ satisfy the
assumptions made in \textit{Step 1} (which holds for large $\nu$)
then the number of the irreducible components of $X$ does not
exceed the number of the irreducible components of $X_{\nu}.$
Therefore by Theorem \ref{manusc} for every irreducible component
$Y$ of $X$ there is a sequence of purely $n$-dimensional analytic
sets $\{Y_{\nu}\}$ converging to $Y$ locally uniformly such that
$Y_{\nu}\subset X_{\nu}$ and $s(\pi|_{Y})=s(\pi|_{Y_{\nu}})$ for
almost all $\nu.$ Consequently, by Lemma~\ref{eqconv},
$\{Y_{\nu}\}$ converges to $Y$ in the sense of holomorphic chains
as required.

To do this we need the following claim. Let $F\subset(U_1\times
B\times\mathbf{C}^k)\cap X$ be the graph of a holomorphic mapping
defined on $U_1\times B.$ (Note that, by (1) and (3), $(U_1\times
B\times\mathbf{C}^k)\cap X$ is the union of such graphs.) Put
$\Sigma=\Sigma(X).$
\begin{claim}\label{claimcom}\em Let
$\tilde{X}=\bigcup_{a\in(U_1\setminus\Sigma)}X_a,$ where for every
$a\in(U_1\setminus\Sigma)$ by $X_a$ we denote the irreducible
component of $(\{a\}\times U_2\times\mathbf{C}^k)\cap X$
containing $(\{a\}\times B\times\mathbf{C}^k)\cap F.$ Then
$\tilde{X}$ is an analytic subset of $((U_1\setminus\Sigma)\times
U_2\times\mathbf{C}^k).$
\end{claim}
\textit{Proof of Claim \ref{claimcom}.} It is sufficient to check
that for every $a\in U_1\setminus\Sigma$ there is a ball
$B'\subset U_1\setminus\Sigma$ centered at $a$ such that
$(B'\times U_2\times\mathbf{C}^k)\cap\tilde{X}$ is an analytic
subset of $B'\times U_2\times\mathbf{C}^k.$

Fix $a_0\in U_1\setminus\Sigma$ and take a ball $B'\subset
U_1\setminus\Sigma$ centered at $a_0.$ We check that $(B'\times
U_2\times\mathbf{C}^k)\cap\tilde{X}$ equals the irreducible
component (denoted by $X'$) of $(B'\times
U_2\times\mathbf{C}^k)\cap X$ containing $(B'\times
B\times\mathbf{C}^k)\cap F.$ First note that $(B'\times
U_2\times\mathbf{C}^k)\cap\tilde{X}\subseteq X'$ (an immediate
consequence of the fact that for every $a\in B'$ the analytic set
$(\{a\}\times U_2\times\mathbf{C}^k)\cap X'$ contains
$(\{a\}\times B\times\mathbf{C}^k)\cap F$ so it must contain $X_a$
as well).

For the converse inclusion, suppose for a moment that
$X'\nsubseteq(B'\times U_2\times\mathbf{C}^k)\cap\tilde{X}.$ Then
there is $(a,b)\in(B'\times B)$ such that the number of points in
$(\{(a,b)\}\times\mathbf{C}^k)\cap X'$ is strictly greater than
the number of points in
$(\{(a,b)\}\times\mathbf{C}^k)\cap\tilde{X}.$ Since $X'$ is
irreducible, there is an arc $$\gamma:[0,1]\rightarrow(((B'\times
U_2)\setminus\pi(\mathcal{S}(X)))\times\mathbf{C}^k)\cap X'$$
connecting all the points in $(\{(a,b)\}\times\mathbf{C}^k)\cap
X'.$

It is easy to see (at least when $B'$ is small which we may
assume) that there is a homeomorphic deformation $H:(B'\times
U_2)\rightarrow(B'\times U_2),$ such that $H(\{a'\}\times
U_2)\subset\{a'\}\times U_2$ for every $a'\in B',$ after which the
set $\pi(\mathcal{S}(X))\cap(B'\times U_2)$ becomes the union of
graphs of constant functions defined on $B'.$ Then the arc
$\tilde{H}\circ\gamma,$ where
$\tilde{H}=(H,id_{\mathbf{C}^k}):B'\times
U_2\times\mathbf{C}^k\rightarrow B'\times U_2\times\mathbf{C}^k,$
can be deformed by shifting along
$\tilde{E}=H(\pi(\mathcal{S}(X))\cap(B'\times U_2))$ to the arc
$$\tau:[0,1]\rightarrow(((\{a\}\times U_2)\setminus\tilde{E})\times\mathbf{C}^k)\cap
\tilde{H}(X').$$ Consequently, $\tilde{H}^{-1}\circ\tau$ is an arc
connecting all the points of $(\{(a,b)\}\times\mathbf{C}^k)\cap
X'$ whose image is contained in $(((\{a\}\times U_2)\setminus
\pi(\mathcal{S}(X)))\times\mathbf{C}^k)\cap X'.$ This means that
$(\{a\}\times U_2\times\mathbf{C}^k)\cap X'$ is irreducible hence
$$(\{a\}\times U_2\times\mathbf{C}^k)\cap\tilde{X}=(\{a\}\times
U_2\times\mathbf{C}^k)\cap X'$$ because $(\{a\}\times
B\times\mathbf{C}^k)\cap F$ is contained in both sets of the
latter equation. On the other hand, these sets have different
number of points in generic fibers over $\{a\}\times U_2,$ a
contradiction.

Thus we have checked that $(B'\times
U_2\times\mathbf{C}^k)\cap\tilde{X}=X'$ which implies the
analyticity and the proof is complete.\qed\vspace*{2mm}

Let us return to the proof of Theorem \ref{component}. For every
irreducible component $Y$ of $X$ select one graph $F_Y$ of a
holomorphic mapping defined on $U_1\times B$ such that
$F_{Y}\subset Y.$ Then, by \textit{Step 1}, there is the uniquely
determined graph $F_{Y,\nu}$ of the holomorphic mapping defined on
$U_1\times B$ such that $F_{Y,\nu}\subset X_{\nu}$ and such that
$$dist((\{x\}\times\mathbf{C}^k)\cap F_Y,(\{x\}\times\mathbf{C}^k)\cap F_{Y,\nu})<\delta$$
for every $x\in U_1\times B$ (for $\delta$ picked in \textit{Step
1}).

Now put $\Sigma_{\nu}=\Sigma(X_{\nu})$ and define
$\tilde{Y}_{\nu}=\bigcup_{a\in U_1\setminus\Sigma_{\nu}}Y_{a,\nu}$
where $Y_{a,\nu}$ is the irreducible component of $(\{a\}\times
U_2\times\mathbf{C}^k)\cap X_{\nu}$ containing $(\{a\}\times
B\times\mathbf{C}^k)\cap F_{Y,\nu}.$ By Claim~\ref{claimcom},
applied to $X_{\nu},$ $\tilde{Y}_{\nu}$ is an analytic subset of
$(U_1\setminus\Sigma_{\nu})\times U_2\times\mathbf{C}^k$ (clearly,
of pure dimension $n$ with proper projection onto
$(U_1\setminus\Sigma_{\nu})\times U_2$). Let $Y_{\nu}$ be the
closure of $\tilde{Y}_{\nu}$ in $U\times\mathbf{C}^k.$ Note that
$Y_{\nu}$ is analytic as $\Sigma_{\nu}\times U_2$ is a nowhere
dense analytic subset of $U$ and $\tilde{Y}_{\nu}\subset X_{\nu}$
and $X_{\nu}$ is with proper projection onto $U$ (so its fibers
over $U$ are locally bounded at every $x\in\Sigma_{\nu}\times
U_2$). It is easy to see that $Y_{\nu}$ is an irreducible
component of $X_{\nu}$ (otherwise $(\{a\}\times
U_2\times\mathbf{C}^k)\cap Y_{\nu}$ would be reducible for some
$a\in U_1\setminus\Sigma_{\nu}$, a contradiction with the
definition of $\tilde{Y}_{\nu}$).

We show that the mapping which assigns to every irreducible
component $Y$ of $X$ the set $Y_{\nu}$ described above is
injective, which completes the proof. To do this it is sufficient
to check the following two facts for
every fixed irreducible component $Y$ of $X$:\vspace*{2mm}\\
(a) $s(\pi|_{Y})=s(\pi|_{Y_{\nu}}),$\\
(b) for $x_0\in U_1\times B$ fixed in \textit{Step 1} and for
every $(x_0,v)\in Y$ there is $(x_0,v_{\nu})\in
Y_{\nu}$\linebreak\hspace*{5.6mm}such that
$||v-v_{\nu}||<\delta.$\vspace*{2mm}\\
Clearly, for every irreducible component $Y_{\nu}$ of $X_{\nu}$
there is at most one $Y$ satisfying (a) and (b).

Let us handle (a). Take $x'_0$ picked in \textit{Step 1}. We may
assume that $x'_0\notin \Sigma_{\nu}$ (otherwise it can be
replaced by a point from an arbitrarily small neighborhood of
$x'_0$ in $U_1\setminus(\Sigma\cup\Sigma_{\nu})$ satisfying all
the conditions specified for $x'_0$ in \textit{Step 1}).

Let $C_Y,C_{Y,\nu}$ denote the irreducible components of
$(\{x'_0\}\times U_2\times\mathbf{C}^k)\cap Y,$ $(\{x'_0\}\times
U_2\times\mathbf{C}^k)\cap X_{\nu}$ respectively containing
$(\{x'_0\}\times B\times\mathbf{C}^k)\cap F_Y,$ $(\{x'_0\}\times
B\times\mathbf{C}^k)\cap F_{Y,\nu}$ respectively.
\begin{claim}\label{cargen}\em The cardinalities of the generic fibers in $C_Y$ and
$C_{Y,\nu}$ over $\{x'_0\}\times U_2$ are equal.
\end{claim}
\textit{Proof of Claim \ref{cargen}.} Observe that the mapping
$$G:(\{x'_0\}\times(U_2\setminus\tilde{B})\times\mathbf{C}^k)\cap Y\rightarrow(\{x'_0\}\times(U_2\setminus\tilde{B})\times\mathbf{C}^k)\cap X_{\nu},$$
where $\tilde{B}=B_1\cup\ldots\cup B_{\hat{k}},$ given by
$G(x'_0,y,u)=(x'_0,y,v),$ where $v$ is the unique vector in
$\mathbf{C}^k$ such that $||u-v||<\delta,$ is a biholomorphism
onto its image. Moreover, since, by the choice of $x'_0,$ the
interior of every $B_i$ contains precisely one element from
$\pi(\mathcal{S}(X))$ and one from $\pi(\mathcal{S}(X_{\nu}))$ and
the intersection of every pair of distinct $B_i$'s is empty, we
may apply the following claim.
\begin{claim}\label{curve}\em Let $D\subset\mathbf{C}$ be a domain,
$\rho:D\times\mathbf{C}^k\rightarrow D$ be the natural projection
and let $E\subset D\times\mathbf{C}^k$ be an irreducible analytic
curve such that $\rho|_{E}$ is proper. Finally, let
$K_1,\ldots,K_s\subset D$ be compact balls, $K_i\cap
K_j=\emptyset$ for $i\neq j$ and let $\sharp
K_i\cap\rho(\mathcal{S}(E))=1$ for every $i=1,\ldots,s.$ Then
$E\cap((D\setminus\bigcup_{i=1}^{s}K_i)\times\mathbf{C}^k)$ is
irreducible.
\end{claim}
\textit{Proof of Claim \ref{curve}.} Put
$$\tilde{E}=E\cap((D\setminus\bigcup_{i=1}^{s}K_i)\times\mathbf{C}^k).$$
It is sufficient to show that for every
$a,b\in\tilde{E}\setminus(\rho(\mathcal{S}(E))\times\mathbf{C}^k)$
there is an arc
$\tau:[0,1]\rightarrow\tilde{E}\setminus(\rho(\mathcal{S}(E))\times\mathbf{C}^k)$
such that $\tau(0)=a,\tau(1)=b.$

Fix
$a,b\in\tilde{E}\setminus(\rho(\mathcal{S}(E))\times\mathbf{C}^k).$
Let $\hat{K}_1,\ldots,\hat{K}_s$ be compact balls in $D,$
$K_i\varsubsetneq\hat{K}_i$ for $i=1,\ldots,s,$ satisfying the
hypotheses of the claim with
$a,b\notin\bigcup_{i=1}^s\hat{K}_i\times\mathbf{C}^k.$
Irreducibility of $E$ implies that there is an arc
$\gamma:[0,1]\rightarrow
E\setminus(\rho(\mathcal{S}(E))\times\mathbf{C}^k)$ such that
$\gamma(0)=a,$ $\gamma(1)=b.$ Then $\tau$ can be obtained from
$\gamma$ as follows. For every $t\in [0,1]$ such that
$\gamma(t)\notin(\bigcup_{i=1}^s\hat{K}_i\times\mathbf{C}^k)$ put
$\tau(t)=\gamma(t).$

Take $t\in [0,1]$ such that
$\gamma(t)\in\hat{K}_i\times\mathbf{C}^k$ for some $i.$ By the
hypothesis $\rho(\mathcal{S}(E))\cap\hat{K}_i=\{g\}$ for some
$g\in D.$ Let $\sigma$ be the segment passing through
$\rho(\gamma(t)),$ connecting $g$ with the uniquely determined
$h(t)\in\partial\hat{K}_i.$ Then $(\sigma\times\mathbf{C}^k)\cap
E$ is the union of graphs of continuous functions defined on
$\sigma.$ Let $f\subset(\sigma\times\mathbf{C}^k)$ be the one of
these graphs for which $\gamma(t)\in f.$ Set
$\tau(t)=(h(t),f(h(t))).$

It is easy to check that $\tau$ defined as above is an arc
connecting $a$ and $b$ in
$\tilde{E}\setminus(\rho(\mathcal{S}(E))\times\mathbf{C}^k).$\qed\vspace*{2mm}\\
\textit{Proof of Claim \ref{cargen} (end).} We show that
$$G((\{x'_0\}\times(U_2\setminus\tilde{B})\times\mathbf{C}^k)\cap C_Y)=(\{x'_0\}\times(U_2\setminus\tilde{B})\times\mathbf{C}^k)\cap C_{Y,\nu},$$
which clearly implies the assertion of the claim. To do this
observe that, by Claim~\ref{curve},
$(\{x'_0\}\times(U_2\setminus\tilde{B})\times\mathbf{C}^k)\cap
C_Y$ is irreducible. Then
$G((\{x'_0\}\times(U_2\setminus\tilde{B})\times\mathbf{C}^k)\cap
C_Y)$ is irreducible as well because $G$ is a biholomorphism onto
its image. This implies that
$$G((\{x'_0\}\times(U_2\setminus\tilde{B})\times\mathbf{C}^k)\cap
C_Y)=(\{x'_0\}\times(U_2\setminus\tilde{B})\times\mathbf{C}^k)\cap
C_{Y,\nu}$$ because
$(\{x'_0\}\times(U_2\setminus\tilde{B})\times\mathbf{C}^k)\cap
C_{Y,\nu}$ is irreducible by Claim \ref{curve} and both sets
contain $(\{x'_0\}\times B\times\mathbf{C}^k)\cap F_{Y,\nu}$
(recall that, by \textit{Step 1,}
$B\cap\bigcup_{i=1}^{\hat{k}}B_i=\emptyset$). \qed\vspace*{2mm}\\
\textit{Proof of Theorem \ref{component} (end).} Now define
$\tilde{Y}=\bigcup_{a\in U_1\setminus\Sigma}Y_a$ where $Y_a$ is
the irreducible component of $(\{a\}\times
U_2\times\mathbf{C}^k)\cap Y$ containing $(\{a\}\times
B\times\mathbf{C}^k)\cap F_Y.$ By Claim \ref{claimcom},
$\tilde{Y}$ is an analytic subset of $(U_1\setminus\Sigma)\times
U_2\times\mathbf{C}^k$ (of pure dimension $n$ with proper
projection onto $(U_1\setminus\Sigma)\times U_2$). Since
$\tilde{Y}$ is contained in $Y\cap((U_1\setminus\Sigma)\times
U_2\times\mathbf{C}^k)$ which is irreducible and $n$-dimensional,
it holds
$$\tilde{Y}=Y\cap((U_1\setminus\Sigma)\times
U_2\times\mathbf{C}^k).$$ The latter fact implies that
$$Y\cap(\{x'_0\}\times U_2\times\mathbf{C}^k)=C_Y$$ and, consequently,
that the cardinalities of the generic fibers in $Y$ over $U$ and
in $C_Y$ over $\{x'_0\}\times U_2$ are equal (because
$\{x'_0\}\times U_2\nsubseteq\pi(\mathcal{S}(X))$).

Finally, observe that (in view of the fact that $\{x'_0\}\times
U_2\nsubseteq\pi(\mathcal{S}(X_{\nu}))$) the cardinalities of the
generic fibers in $Y_{\nu}$ over $U$ and in $C_{Y,\nu}$ over
$\{x'_0\}\times U_2$ are equal because
$$Y_{\nu}\cap(\{x'_0\}\times U_2\times\mathbf{C}^k)=C_{Y,\nu}$$ by
the definition of $Y_{\nu}.$ Thus, in view of Claim \ref{cargen},
$s(\pi|_{Y})=s(\pi|_{Y_{\nu}})$ as required in (a).

Let us turn to (b). Consider the arc $\gamma_{Y}$ defined in
\textit{Step 1.} Note that for every $\gamma_Y(t)\in X,$ where
$t\in [0,1],$ there is precisely one element $(e(t),f(t))\in
X_{\nu}\subset U\times\mathbf{C}^k$ such that
$\pi(\gamma_Y(t))=e(t)$ and the distance between $(e(t),f(t))$ and
$\gamma_Y(t)$ is smaller than $\delta.$ Since
$(\pi(\gamma_Y([0,1])))\cap\pi(\mathcal{S}(X_{\nu}))=\emptyset$
then
$$\tau:[0,1]\ni t\mapsto(e(t),f(t))\in X_{\nu}$$ is an arc whose image
is contained in one irreducible component of $X_{\nu}.$ On the
other hand, there is $t_0$ such that $\gamma_Y(t_0)\in F_Y$ so
$\tau(t_0)\in F_{Y,\nu},$ which implies that the irreducible
component containing $\tau([0,1])$ contains $F_{Y,\nu}$ as well.
Thus this irreducible component must be $Y_{\nu}.$ To complete the
proof observe that for every $(x_0,v)\in Y$ there is $t'\in [0,1]$
such that $(x_0,v)=\gamma_Y(t'),$ $\tau(t')=(x_0,f(t'))\in
Y_{\nu},$ $||v-f(t')||<\delta.$

Thus we have proved the theorem under the extra hypotheses (1),
(2) and (3) formulated at the beginning.\vspace*{2mm}\\
\textit{Step 3.} Let us show that (1), (2) and (3) need not be
assumed. Let $\{X_{\nu}\}$ be a sequence of analytic sets
satisfying the hypotheses of Theorem \ref{component} and let $Y$
be an analytic subset of $U\times\mathbf{C}^k$ of pure dimension
$n$ with $Y\subset X.$ Fix an open set $\tilde{U}\subset\subset
U.$

Cover $\overline{\tilde{U}}$ by a finite number of domains
$E_1,\ldots,E_s\subset U,$ $E_i\subset\subset\tilde{E}_i\subset U$
such that for every $i\in\{1,\ldots,s\}$ there is a polynomial
automorphism $L_i:\mathbf{C}^n\rightarrow\mathbf{C}^n$ with the
following property. The conditions (1), (2) and (3) are satisfied
with $U, X$ replaced by
$L_i(\tilde{E}_i),\beta_{L_i}(X\cap(\tilde{E}_i\times\mathbf{C}^k))$
respectively, where
$\beta_{L_i}:\mathbf{C}^n\times\mathbf{C}^k\rightarrow\mathbf{C}^n\times\mathbf{C}^k$
is given by the formula $\beta_{L_i}(x,z)=(L_i(x),z).$

By \textit{Step 2} for every $i\in\{1,\ldots,s\}$ there is a
sequence $\{Y_{i,\nu}\}$ of purely\linebreak $n$-~dimensional
analytic subsets of $E_i\times\mathbf{C}^k,$
$Y_{i,\nu}\subset(E_i\times\mathbf{C}^k)\cap X_{\nu},$ converging
to $(E_i\times\mathbf{C}^k)\cap Y$ in the sense of chains. Let us
check that
$Y_{\nu}=(\bigcup^s_{i=1}Y_{i,\nu})\cap(\tilde{U}\times\mathbf{C}^k)$
is, for large $\nu,$ an analytic subset of
$\tilde{U}\times\mathbf{C}^k.$ One easily observes that this is
the case as the convergence in the sense of chains of
$\{X_{\nu}\}$ to $X$ imply that for almost all $\nu$ and for every
$i,j$ it holds
$$Y_{i,\nu}\cap((E_i\cap
E_j)\times\mathbf{C}^k)=Y_{j,\nu}\cap((E_i\cap
E_j)\times\mathbf{C}^k).$$ The latter equation also implies that
$\{Y_{\nu}\}$ converges to $Y\cap(\tilde{U}\times\mathbf{C}^k)$ in
the sense of chains. Thus the proof of Theorem \ref{component} is
complete.\qed
\subsection{Approximation of mappings}\label{mainproof}
%
First let us note that in Theorem \ref{main} (and in its
generalizations) the space $\mathbf{C}^n$ containing $U$ may be
replaced by an affine algebraic variety. In fact, in the global
version of the approximation theorem (see \cite{Lem}, Theorem~3.2)
the domain of the approximated mapping is admitted to have
singularities. Since this case is reduced to the one where the
mapping is defined on an open subset of some $\mathbf{C}^n$ and
the reduction is of purely analytic geometric nature, we assume
here that $U$ is an open subset of $\mathbf{C}^n.$

Our aim is to give a direct geometric proof of Theorem \ref{main},
or more precisely, its semi-global version (Theorem
\ref{semiglobmain}). The proof is organized as follows. First
using Theorem \ref{component} we prove in Section \ref{globproven}
the following

\begin{proposition}\label{proglo}Let $U$ be a domain in $\mathbf{C}^n$ and
let $f:U\rightarrow\mathbf{C}^k$ be a holomorphic mapping that
satisfies a system of equations $Q(x,f(x))=0$ for $x\in U.$ Here
$Q$ is a Nash mapping from a neighborhood in
$\mathbf{C}^n\times\mathbf{C}^k$ of the graph of $f$ into some
$\mathbf{C}^q.$ Then there is $R\in\mathbf{C}[x,z_1,\ldots,z_k]$
with $R(x,f(x))$ not identically zero such that the following
holds. If for some open $\tilde{U}, U_0\subset U$ with
$U_0\subset\subset\tilde{U}$ there is a sequence
$\{g^{\nu}:\tilde{U}\rightarrow\mathbf{C}^k\}$ of Nash mappings
converging locally uniformly to $f|_{\tilde{U}}$ such that
$\{\{x\in\tilde{U}: R(x,g^{\nu}(x))=0\}\}$ converges to
$\{x\in\tilde{U}: R(x,f(x))=0\}$ in the sense of chains then there
is a sequence $\{f^{\nu}:U_0\rightarrow\mathbf{C}^k\}$ of Nash
mappings converging uniformly to $f|_{U_0}$ such that
$Q(x,f^{\nu}(x))=0$ for $x\in U_0,$ $\nu\in\mathbf{N}.$
\end{proposition}

Next in Section \ref{construction} for any holomorphic mapping
$f:U\rightarrow\mathbf{C}^k,$ $f=f(x),$ and any
$R\in\mathbf{C}[x,z_1,\ldots,z_k]$ such that $R(x,f(x))$ is not
identically zero, we construct a family $\mathcal{U}_{f,R}$ of
open subsets of $U$ such that the following holds.

\begin{proposition}\label{propoly}
Let $f:U\rightarrow\mathbf{C}^k$ be a holomorphic mapping, where
$U$ is a domain in $\mathbf{C}^n,$ let
$R\in\mathbf{C}[x,z_1,\ldots,z_k]$ be such that $R(x,f(x))$ is not
identically zero on $U$ and let $U_0\in\mathcal{U}_{f,R}.$ Then
there are an open $\tilde{U}\subset\mathbf{C}^n$ with
$U_0\subset\subset\tilde{U}$ and a sequence
$\{f^{\nu}:\tilde{U}\rightarrow\mathbf{C}^k\}$ of Nash mappings
converging uniformly to $f|_{\tilde{U}}$ such that
$\{\{x\in\tilde{U}: R(x,f^{\nu}(x))=0\}\}$ converges to
$\{x\in\tilde{U}: R(x,f(x))=0\}$ in the sense of chains.
\end{proposition}

Proposition \ref{propoly} is proved in Subsection
\ref{globperties}. One of the main results of this paper is the
following semi-global version of the approximation theorem.

\begin{theorem}\label{semiglobmain}Let $U$ be a domain in $\mathbf{C}^n$ and
let $f:U\rightarrow\mathbf{C}^k$ be a holomorphic mapping that
satisfies a system of equations $Q(x,f(x))=0$ for $x\in U.$ Here
$Q$ is a Nash mapping from a neighborhood in
$\mathbf{C}^n\times\mathbf{C}^k$ of the graph of $f$ into some
$\mathbf{C}^q.$ Let $R$ be any polynomial obtained by applying
Proposition \ref{proglo} with $f, Q.$ Then for every
$U_0\in\mathcal{U}_{f,R}$ there is a sequence
$\{f^{\nu}:U_0\rightarrow\mathbf{C}^k\}$ of Nash mappings
converging uniformly to $f|_{U_0}$ such that $Q(x,f^{\nu}(x))=0$
for $x\in U_0$ and $\nu\in\mathbf{N}.$
\end{theorem}
\proof In view of the fact that $R$ satisfies the assertion of
Proposition \ref{proglo}, it is sufficient to apply
Proposition~\ref{propoly}.\qed\vspace*{2mm}

In order to characterize those $U_0$ for which the presented
methods work we need an insight into $\mathcal{U}_{f,R}.$ Here let
us just mention two special cases which directly follow by
Subsection \ref{construction}

\begin{corollary}(a) Let $U$ be a domain in $\mathbf{C}.$ Then
for every holomorphic mapping $f:U\rightarrow\mathbf{C}^k$ and
every $R\in\mathbf{C}[x,z_1,\ldots,z_k],$ $R(x,f(x))$ not
identically zero, the family $\mathcal{U}_{f,R}$ contains every
open set $U_0$ for which there is a Runge domain $\tilde{U}$ with
$U_0\subset\subset\tilde{U}\subset U.$ Consequently, if $f$
depends on one variable we have the global version of
Theorem~\ref{main}, first proved by van
den Dries in \cite{vD}.\vspace*{2mm}\\
(b) Let $U$ be a domain in $\mathbf{C}^n.$ Then for every
holomorphic mapping $f:U\rightarrow\mathbf{C}^k,$ every
$R\in\mathbf{C}[x,z_1,\ldots,z_k],$ $R(x,f(x))$ not identically
zero, and $x_0\in U$ there is an open neighborhood $\tilde{U}$ of
$x_0$ such that $\tilde{U}\in\mathcal{U}_{f,R},$ which implies
Theorem \ref{main}.
\end{corollary}

\subsubsection{Proof of Proposition \ref{proglo}}\label{globproven}
Proposition \ref{proglo} is a consequence of the following
\begin{proposition}\label{claglo}Let $U, V$ be a domain in $\mathbf{C}^n$ and an algebraic subvariety of
$\mathbf{C}^{\hat{m}}$ respectively. Let $F:U\rightarrow V$ be a
holomorphic mapping. Then there is a polynomial $R$ in $\hat{m}$
variables with $R\circ F$ not identically zero such that the
following holds. If for some open $\tilde{U}, U_0\subset U$ with
$U_0\subset\subset\tilde{U}$ there is a sequence
$\{G^{\nu}:\tilde{U}\rightarrow\mathbf{C}^{\hat{m}}\}$ of Nash
mappings converging locally uniformly to $F|_{\tilde{U}}$ such
that $\{\{x\in\tilde{U}: (R\circ G^{\nu})(x)=0\}\}$ converges to
$\{x\in\tilde{U}: (R\circ F)(x)=0\}$ in the sense of chains then
there is a sequence $\{F^{\nu}:U_0\rightarrow V\}$ of Nash
mappings converging uniformly to $F|_{U_0}.$
\end{proposition}

First let us check that Proposition \ref{claglo} implies
Proposition \ref{proglo}. Let $f:U\rightarrow\mathbf{C}^k$ be the
holomorphic mapping from Proposition \ref{proglo}. Put
$F(x)=(x,f(x))$ and $\hat{m}=n+k.$ Let $V$ be the intersection of
all algebraic subvarieties of $\mathbf{C}^{\hat{m}}$ containing
$F(U).$ Then by Proposition \ref{claglo} there is
$R\in\mathbf{C}[x_1,\ldots,x_n,z_1,\ldots,z_k]$ satisfying the
assertion of this proposition. Next fix $U_0,\tilde{U}$ as in
Proposition~\ref{proglo}, assume (without loss of generality) that
$U_0$ is connected and take an open connected $U_1$ with
$U_0\subset\subset U_1\subset\subset\tilde{U}.$ Now let
$g^{\nu}:\tilde{U}\rightarrow\mathbf{C}^k$ be a sequence of Nash
mappings converging locally uniformly to $f|_{\tilde{U}}$ such
that $\{\{x\in\tilde{U}: R(x,g^{\nu}(x))=0\}\}$ converges to
$\{x\in\tilde{U}: R(x,f(x))=0\}$ in the sense of chains. Set
$G^{\nu}(x)=(x,g^{\nu}(x)).$ Then by Proposition \ref{claglo}
there is a sequence $\{F^{\nu}:U_1\rightarrow V\}$ of Nash
mappings converging uniformly to $F|_{U_1}.$

We need to show that the first $n$ components of $F^{\nu}$ may be
assumed to constitute the identity and that $Q\circ F^{\nu}=0$ for
sufficiently large $\nu$. To this end denote
$Y=\{(x,v)\in\hat{U}\subset\mathbf{C}^n\times\mathbf{C}^k:Q(x,v)=0\},$
where $\hat{U}$ is the domain of $Q.$ Clearly, we may assume that
$F^{\nu}(U_1)\subset\hat{U}$ for almost all $\nu.$ Next observe
that $F^{\nu}(U_1)\subset Y$ for almost all $\nu.$ Indeed, take
$\hat{z}\in F(U_1)\cap Reg(V)$ (the intersection is non-empty as
$F(U_1)\subset Sing(V)$ implies, by the connectedness of $U,$ that
$F(U)\subset Sing(V)\varsubsetneq V $). Let $B$ be an open
neighborhood of $\hat{z}$ in $\mathbf{C}^n\times\mathbf{C}^k$ such
that $B\cap V$ is a connected manifold and let $U_2$ be a
non-empty open subset of $U_1$ such that $F(U_2),
F^{\nu}(U_2)\subset B$ for almost all $\nu.$ Then $B\cap V\subset
Y$ (otherwise $F(U_2)\subset\tilde{V}$ where $\tilde{V}$ is an
algebraic subvariety of $\mathbf{C}^n\times\mathbf{C}^k$ with
$dim(\tilde{V})<dim(V)$). This implies that $F^{\nu}(U_2)\subset
Y$ for almost all $\nu$ hence $F^{\nu}(U_1)\subset Y$ because
$U_1$ is connected.

Let $\tilde{F}^{\nu}:U_1\rightarrow\mathbf{C}^n,$ for
$\nu\in\mathbf{N},$ be the mapping whose components are the first
$n$ components of $F^{\nu}.$ Since $\{\tilde{F}^{\nu}\}$ converges
uniformly to the identity on $U_1$ and $U_0\subset\subset U_1$
there is a sequence $H^{\nu}:U_0\rightarrow U_{1}$ of Nash
mappings such that $\tilde{F}^{\nu}\circ H^{\nu}=id_{U_0}$ if
$\nu$ is large enough. Consequently, $F^{\nu}\circ
H^{\nu}(x)=(x,f^{\nu}(x))$ for $x\in U_0$ and
$\{f^{\nu}:U_0\rightarrow\mathbf{C}^k\}$ satisfies the assertion
of Proposition \ref{proglo}.\vspace*{2mm}\\
\textit{Proof of Proposition \ref{claglo}.} First observe that we
may assume $F(U)\nsubseteq Sing(V)$ as otherwise $V$ may be
replaced by $Sing(V).$ Next, since $U$ is connected, $F(U)$ is
contained in one irreducible component of $V$ so we may assume
that $V$ is of pure dimension, say $m.$

We may also assume that
$V\subset\mathbf{C}^{\hat{m}}\thickapprox\mathbf{C}^m\times\mathbf{C}^{s}$
is with proper projection onto $\mathbf{C}^m.$ Indeed, there is a
$\mathbf{C}$-linear isomorphism
$J:\mathbf{C}^{m+s}\rightarrow\mathbf{C}^{m+s}$ such that $J(V)$
is with proper projection onto $\mathbf{C}^m.$ Thus if there
exists a sequence $H^{\nu}:U_0\rightarrow J(V)$ of Nash mappings
converging to $J\circ F|_{U_0}$ then the sequence $\{J^{-1}\circ
H^{\nu}\}$ satisfies the assertion of the proposition.

To complete the preparations, by
$\rho:\mathbf{C}^m\times\mathbf{C}^s\rightarrow\mathbf{C}^m,$
$\tilde{\rho}:\mathbf{C}^m\times\mathbf{C}\rightarrow\mathbf{C}^m$
denote the natural projections. Passing to the image of $V$ by a
linear isomorphism arbitrarily close to the identity, if
necessary, we assume (in view of $F(U)\nsubseteq Sing(V)$) that
$\rho(F(U))\nsubseteq\rho(\mathcal{S}(V)).$ Now the polynomial $R$
is constructed as follows. Any $\mathbf{C}$-linear form
$L:\mathbf{C}^s\rightarrow\mathbf{C}$ determines the mapping
$\Phi_{L}:\mathbf{C}^m\times\mathbf{C}^s\rightarrow\mathbf{C}^m\times\mathbf{C}$
by $\Phi_L(u,v)=(u,L(v)).$ Since $V$ is an algebraic subset of
$\mathbf{C}^m\times\mathbf{C}^s$ with proper projection onto
$\mathbf{C}^m$ then $\Phi_L(V)$ is an algebraic subset of
$\mathbf{C}^m\times\mathbf{C}$ also with proper projection onto
$\mathbf{C}^m$ for every form $L.$ Take a form $L$ such that the
fibers of the projections of $\Phi_L(V)$ and $V$ onto
$\mathbf{C}^m$ have generically the same cardinality and
$\rho(F(U))\nsubseteq\tilde{\rho}(\mathcal{S}(\Phi_L(V)))$. The
set $\Phi_L(V)$ is described by the unitary polynomial in one
variable (corresponding to the last coordinate of
$\mathbf{C}^m\times\mathbf{C}$) whose coefficients are polynomials
in $m$ variables and whose discriminant is non-zero. The
discriminant, denoted by $R,$ is the polynomial we look for. In
fact, after the preparations, $R$ depends only on $m\leq\hat{m}$
variables (the last $s=\hat{m}-m$ variables are dummy).

Let us show that $R$ indeed has all the required properties. First
$R\circ F$ is not identically zero as
$\rho(F(U))\nsubseteq\tilde{\rho}(\mathcal{S}(\Phi_L(V))).$ Next
take $U_0,$ $\tilde{U}$ and $G^{\nu}$ as in Proposition
\ref{claglo}. We need the following notation. For any holomorphic
mapping $H:E\rightarrow\mathbf{C}^m,$ where $E,E'$ are open
subsets of $\mathbf{C}^n,$ $E'\subset E,$ and any algebraic
subvariety $X$ of $\mathbf{C}^m\times\mathbf{C}^s$ denote
$$\mathcal{V}(X,E',H)=\{(x,v)\in E'\times\mathbf{C}^s:(H(x),v)\in X\}.$$
The mappings $F,G^{\nu}$ are of the form $F=(\tilde{F},\hat{F}),$
$G^{\nu}=(\tilde{G}^{\nu},\hat{G}^{\nu}),$ for some holomorphic
$\tilde{F}: U\rightarrow\mathbf{C}^m,$
$\tilde{G}^{\nu}:\tilde{U}\rightarrow\mathbf{C}^m,$ $\hat{F}:
U\rightarrow\mathbf{C}^s,$
$\hat{G}^{\nu}:\tilde{U}\rightarrow\mathbf{C}^s.$

In order to prove Proposition \ref{claglo} it is sufficient to
prove the following
\begin{claim}\label{clatwo}For every irreducible component $Y$ of $\mathcal{V}(V,U_0,\tilde{F})$
there is a sequence $\{Y_{\nu}\}$ of analytic subsets of
$U_0\times\mathbf{C}^s$ converging to $Y$ in the sense of chains
such that $Y_{\nu}$ is an irreducible component of
$\mathcal{V}(V,U_0,\tilde{G}^{\nu})$ for every $\nu.$
\end{claim}
For the notion of the convergence of holomorphic chains see
Section \ref{holchai}. Let us check that Proposition \ref{claglo}
indeed follows by Claim \ref{clatwo}. To this end note that
$graph(\hat{F})\subset\mathcal{V}(V,U_0,\tilde{F}).$ Since there
is a sequence $\{Y_{\nu}\}$ of analytic sets converging to
$graph(\hat{F})$ in the sense of chains then
$Y_{\nu}=graph(H^{\nu})$ for almost all $\nu,$ where
$H^{\nu}:U_0\rightarrow\mathbf{C}^s$ is a holomorphic mapping. In
fact, since $Y_{\nu}$ is an irreducible component of
$\mathcal{V}(V,U_0,\tilde{G}^{\nu})$ which is a Nash set then
$Y_{\nu}$ is a Nash set as well. Consequently $H^{\nu}$ is a Nash
mapping. Obviously, $H^{\nu}$ converges to $\hat{F}$  so
$(\tilde{G}^{\nu},H^{\nu}):U_0\rightarrow V$ converges to
$F|_{U_0}$ as required.

Let us turn to the proof of Claim \ref{clatwo}. First let us show
that it is sufficient to prove this claim in the case where $V$ is
replaced by $\Phi_L(V)$ where
$L:\mathbf{C}^s\rightarrow\mathbf{C}$ is the linear form which has
been used to define $R.$ By analogy to the definition of $\Phi_L$
put $\Psi_L(x,v)=(x,L(v))$ for any $x\in\mathbf{C}^n,
v\in\mathbf{C}^s.$ Let
$\tilde{\pi}:\mathbf{C}^n\times\mathbf{C}\rightarrow\mathbf{C}^n,$
$\pi:\mathbf{C}^n\times\mathbf{C}^s\rightarrow\mathbf{C}^n$ denote
the natural projections. We need the following obvious
\begin{remark}\label{funnyobvious}\em
Let $Z\subset E \times\mathbf{C}^s$ be an analytic subset of pure
dimension $n$ with proper projection onto a domain
$E\subset\mathbf{C}^n$ such that
$s(\pi|_Z)=s(\tilde{\pi}|_{\Psi_{L}(Z)}).$ Then for every
irreducible analytic component $\Sigma$ of $\Psi_{L}(Z)$ there
exists an irreducible analytic component $\Gamma$ of $Z$ such that
$\Psi_{L}(\Gamma)=\Sigma$ and
$s(\pi|_{\Gamma})=s(\tilde{\pi}|_{\Sigma}).$
\end{remark}

Assume that Claim \ref{clatwo} holds with $\Phi_{L}(V)$ taken
instead of $V$ ($s=1$). We check that it also holds with $V.$
First observe that
$\mathcal{V}(\Phi_{L}(V),U_0,\tilde{F})=\Psi_{L}(\mathcal{V}(V,U_0,\tilde{F}))$
and
$\mathcal{V}(\Phi_{L}(V),U_0,\tilde{G}^{\nu})=\Psi_{L}(\mathcal{V}(V,U_0,\tilde{G}^{\nu}))$
for $\nu\in\mathbf{N}$ and fix an irreducible component $Y$ of
$\mathcal{V}(V,U_0,\tilde{F}).$ Then there are irreducible
components $\Theta_{\nu}$ of
$\Psi_{L}(\mathcal{V}(V,U_0,\tilde{G}^{\nu})),$ for
$\nu\in\mathbf{N},$ such that $\{\Theta_{\nu}\}$ converges to
$\Psi_{L}(Y)$ in the sense of holomorphic chains.

Next note that the fact that
$\tilde{F}(U_0)\nsubseteq\tilde{\rho}(\mathcal{S}(\Phi_L(V)))$ and
the way $L$ has been chosen imply that the cardinalities of the
generic fibers in $\Psi_{L}(\mathcal{V}(V,U_0,\tilde{F})),$
$\mathcal{V}(V,U_0,\tilde{F}),$
$\Psi_{L}(\mathcal{V}(V,U_0,\tilde{G}^{\nu}))$ and in
$\mathcal{V}(V,U_0,\tilde{G}^{\nu})$ over $U_0$ are equal for
large $\nu.$ Therefore, by Remark~\ref{funnyobvious}, for almost
all $\nu$ there is an irreducible component $Y_{\nu}$ of
$\mathcal{V}(V,U_0,\tilde{G}^{\nu})$ such that
$\Psi_{L}(Y_{\nu})=\Theta_{\nu}$ and
$s(\pi|_{Y_{\nu}})=s(\pi|_{Y}).$ Thus it remains to check, in view
of Lemma \ref{eqconv}, that $\{Y_{\nu}\}$ converges to $Y$ locally
uniformly. Observe that otherwise there would be a subsequence
$\{Y_{\nu_{\mu}}\}$ of $\{Y_{\nu}\}$ converging to a purely
$n$-dimensional analytic set $Z\neq Y.$ But then, by the fact that
$\Psi_{L}$ preserves the cardinality of the generic fiber in
$\mathcal{V}(V,U_0,\tilde{F}),$ it holds
$\Psi_{L}(Z)\neq\Psi_{L}(Y)$ which contradicts the fact that
$\{\Psi_{L}({Y_{\nu}})\}$ converges to $\Psi_{L}(Y).$

Now we turn to the proof of Claim \ref{clatwo} with $V$ replaced
by $\Phi_L(V).$ It holds
$$\Phi_L(V)=\{(y,z)\in\mathbf{C}^m\times\mathbf{C}:P(y,z)=0\}$$ where
$$P(y,z)=z^t+z^{t-1}c_1(y)+\ldots+c_t(y)\in(\mathbf{C}[y])[z]$$ for some $t\in\mathbf{N}.$
We may assume that $P$ treated as a polynomial in $z$ has a
non-zero discriminant. Let us recall that the polynomial $R$ is,
by definition, this discriminant.

To complete the proof put
$$X_{\nu}=\mathcal{V}(\Phi_L(V),\tilde{U},\tilde{G}^{\nu})=\{(x,z)\in
\tilde{U}\times\mathbf{C}:P(\tilde{G}^{\nu}(x),z)=0\}$$ and
$$X=\mathcal{V}(\Phi_L(V),\tilde{U},\tilde{F})=\{(x,z)\in
\tilde{U}\times\mathbf{C}:P(\tilde{F}(x),z)=0\}.$$ Then
$$\tilde{\pi}(\mathcal{S}(X_{\nu}))=\{x\in\tilde{U}: R(\tilde{G}^{\nu}(x))=0\}$$
and
$$\tilde{\pi}(\mathcal{S}(X))=\{x\in\tilde{U}:
R(\tilde{F}(x))=0\},$$ hence by the hypothesis the sequence
$\{\tilde{\pi}(\mathcal{S}(X_{\nu}))\}$ converges to
$\tilde{\pi}(\mathcal{S}(X))$ in the sense of holomorphic chains.
Now it is sufficient to apply Theorem \ref{component} and the
proof of Claim \ref{clatwo} is complete. Consequently, we have
also proved Proposition~\ref{claglo} and Proposition
\ref{proglo}.\qed

\subsubsection{Construction of $\mathcal{U}_{f,R}$}\label{construction}

Put $x=(x_1,\ldots,x_n),$ $x'=(x_1,\ldots,x_{n-1})$ and
$\pi(x)=x'.$ Let $f:U\rightarrow\mathbf{C}^k,$ $f=f(x),$ be a
holomorphic mapping where $U$ is a domain in $\mathbf{C}^n$ and
let $R\in\mathbf{C}[x,z_1,\ldots,z_k]$ be such that $R(x,f(x))$ is
not identically zero. We construct the family $\mathcal{U}_{f,R}.$
The construction is recursive with respect to the number $n$ of
the variables $f$ depends on.

Let $U_0$ be an open subset of $\mathbf{C},$ $n=1.$ Then
$U_0\in\mathcal{U}_{f,R}$ iff $U_0$ is a relatively compact subset
of some open simply connected subset of $U$ (hence in this case
$\mathcal{U}_{f,R}$ depends only on $U$).

Now assume that $U_0$ is an open subset of $\mathbf{C}^n$ for
$n>1.$ Then $U_0\in\mathcal{U}_{f,R}$ iff there is a
biholomorphism $\phi:\check{U}\rightarrow\hat{U}\subset U,$ where
$\hat{U},\check{U}$ are a domain and a Runge domain respectively
in $\mathbf{C}^n$ with $U_0\subset\subset\hat{U},$ and there is a
domain $\check{U}_1\subset\mathbf{C}^{n-1}$ with
$\check{U}\subset\check{U}_1\times\mathbf{C}$ such that the
following hold:\vspace*{2mm}\\
(1) $R(\phi(x),f(\phi(x)))=\tilde{H}(x)W(x),$ $x\in\check{U},$ for
some $\tilde{H}\in\mathcal{O}(\check{U})$ non-vanishing on
$\overline{\phi^{-1}(U_0)}$ and some unitary polynomial
$W\in\mathcal{O}(\check{U}_1)[x_n]$  such
that $W^{-1}(0)\subset\check{U},$\\
(2) $\pi(\phi^{-1}(U_0))\in\mathcal{U}_{g,S}$ for some holomorphic
mapping $g:\check{U}_1\rightarrow \mathbf{C}^s,$ $g=g(x'),$ and
some $S\in\mathbf{C}[x',z_1,\ldots,z_s]$ determined by $f,R,\phi,W$ below.\vspace*{2mm}\\
Given $f, R, \phi, W$ we obtain $g, S$ as follows. Put
$\tilde{f}=f\circ\phi.$ Then $\tilde{f},\phi$ are of the form
$\tilde{f}=(f_1,\ldots,f_k),$ $\phi=(\phi_1,\ldots,\phi_n)$ for
some $f_j,\phi_i\in\mathcal{O}(\check{U})$ for $j=1,\ldots,k,
i=1,\ldots,n.$ By (1) we have:\\
$$f_j(x)=W(x)H_j(x)+r_j(x),$$
$$\phi_i(x)=W(x)\check{H}_i(x)+\check{r}_i(x),$$
for $x\in\check{U},$ where
$r_j(x),\check{r}_i(x)\in\mathcal{O}(\check{U}_1)[x_n]$ satisfy
$deg(r_j), deg(\check{r}_i)<deg(W)$ and
$H_j,\check{H}_i\in\mathcal{O}(\check{U})$ for $j=1,\ldots,k,$
$i=1,\ldots,n.$

Next, there are optimal polynomials (for the definition consult
Section \ref{setswithprop})
$W_1,\ldots,W_{\hat{s}}\in\mathcal{O}(\check{U}_1)[x_n]$ such that
$W=W_1^{k_1}\cdot\ldots\cdot W_{\hat{s}}^{k_{\hat{s}}}$ and
$dim(W_i^{-1}(0)\cap W_j^{-1}(0))<n-1$ for every $i\neq j.$ Put
$d=deg(W).$ For $l=1,\ldots,\hat{s},$ $j=1,\ldots,k,$
$i=1,\ldots,n,$ the polynomials $W_l,r_j,\check{r}_i$ are of the
form:\vspace*{2mm}\\
$W_l(x)=x_n^{p_l}+x_n^{p_l-1}a_{l,1}(x')+\ldots+a_{l,p_l}(x'),$\\
$r_j(x)=x_n^{d-1}b_{j,0}(x')+x_n^{d-2}b_{j,1}(x')+\ldots+b_{j,d-1}(x'),$\\
$\check{r}_i(x)=x_n^{d-1}c_{i,0}(x')+x_n^{d-2}c_{i,1}(x')+\ldots+c_{i,d-1}(x').$\vspace*{2mm}\\
Let $s$ denote the number of all the coefficients of $W_l, r_j,
\check{r}_i$ for all admissible $l, j, i.$ The mapping
$g:\check{U}_1\rightarrow\mathbf{C}^s$ is defined by:
$$g=(A_1,\ldots,A_{\hat{s}},B_{1},\ldots,B_{k},C_{1},\ldots,C_n)$$
where $A_l=(a_{l,1},\ldots,a_{l,p_l}),$
$B_j=(b_{j,0},\ldots,b_{j,d-1}),$ $C_i=(c_{i,0},\ldots,c_{i,d-1})$
again for all admissible $l, j, i.$

Let us turn to determining $S.$ Replacing the holomorphic
coefficients $$a_{l,1},\ldots,a_{l,p_l},b_{j,0},\ldots,b_{j,d-1},
c_{i,0},\ldots,c_{i,d-1}$$ for all $l,$ $j,$ $i$ in
$W_l,r_j,\check{r}_i$ by new variables denoted by the same letters
we obtain polynomials $P_l,w_j,\check{w}_i$ respectively. Put
$P=P_1^{k_1}\cdot\ldots\cdot P_{\hat{s}}^{k_{\hat{s}}}$ and
define:
$$\alpha_j=PS_j+w_j, \beta_i=P\check{S}_i+\check{w}_i$$ for
$j=1,\ldots,k$ and $i=1,\ldots,n,$ where $S_j,\check{S}_i$ are new
variables. Now divide
$R(\beta_1,\ldots,\beta_n,\alpha_1,\ldots,\alpha_k)$ by $P$
(treated as a polynomial in $x_n$ with polynomial coefficients) to
obtain\vspace*{2mm}\\
(*)\hspace*{7mm}$R(\beta_1,\ldots,\beta_n,\alpha_1,\ldots,\alpha_k)=\tilde{W}P+x_n^{d-1}T_1+x_n^{d-2}T_2+\ldots+T_{d},$\vspace*{2mm}\\
where $\tilde{W},T_1,\ldots,T_{d}$ are polynomials such that
$T_1,\ldots,T_{d}$ depend only on the tuple of variables $u,$
where
$$u=(A_1,\ldots,A_{\hat{s}},B_{1},\ldots,B_{k},C_{1},\ldots,C_n)$$
and $A_l=(a_{l,1},\ldots,a_{l,p_l}),$
$B_j=(b_{j,0},\ldots,b_{j,d-1}),$ $C_i=(c_{i,0},\ldots,c_{i,d-1})$
for all admissible $l, j, i.$

Finally, put $T(u)=(T_1(u),\ldots,T_{d}(u))$ and observe that
$T(g(x'))=0$ for $x'\in\check{U}_1.$ By Proposition \ref{proglo}
there is $S\in\mathbf{C}[x',z_1,\ldots,z_s]$ satisfying the
assertion of Proposition \ref{proglo} with $g, T, S$ taken in
place of $f, Q, R.$ Any such $S$ is suitable for our recursive
definition.
\begin{remark}\em
For every holomorphic mapping $f:\mathbf{C}^n\supset
U\rightarrow\mathbf{C}^k,$ polynomial
$R\in\mathbf{C}[x,z_1,\ldots,z_k]$ with $R(x,f(x))$ not
identically zero and every $x_0\in U$ the following holds. There
is a neighborhood $E$ of $x_0$ in $U$ such that $\{x\in E:
R(x,f(x))=0\}$ is either empty or with proper projection onto an
open subset of some affine $n-1$ dimensional subspace of
$\mathbf{C}^{n}.$ This fact, applied recursively in the
construction above, immediately implies that there is an open
neighborhood $\tilde{U}$ of $x_0$ in $U$ such that
$\tilde{U}\in\mathcal{U}_{f,R}.$
\end{remark}

\subsubsection{Proof of Proposition \ref{propoly}}\label{globperties}
The proposition is proved by induction on $n$ (the number of the
variables $f$ depends on). First suppose that $U\subset\mathbf{C}$
(i.e. $n=1$) and fix $f,R$ satisfying the assumptions of the
proposition. Let $U_0$ be any open relatively compact subset of
some open simply connected $\tilde{U}\subset U.$ Then
$$R(x,f(x))=(x-x_0)^{\alpha_0}\cdot\ldots\cdot(x-x_m)^{\alpha_m}g(x)$$
for some $m,\alpha_0,\ldots,\alpha_m\in\mathbf{N},$
$g\in\mathcal{O}(U)$ such that $g(x)\neq 0$ for
$x\in\overline{U_0}.$

Put $W(x)=(x-x_0)^{\alpha_0}\cdot\ldots\cdot(x-x_m)^{\alpha_m}.$
The mapping $f$ is of the form $f=(f_1,\ldots,f_k)$ for some
$f_j\in\mathcal{O}(U),$ $j=1,\ldots,k.$ It holds
$$f_j(x)=W(x)H_j(x)+r_j(x),\mbox{ }x\in U,$$ where $H_j\in\mathcal{O}(U),$ $r_j\in\mathbf{C}[x]$ for $j=1,\ldots,k.$
Now define $f^{\nu}=(f^{\nu}_1,\ldots,f^{\nu}_k)$ on $\tilde{U}$
by $f^{\nu}_j(x)=W(x)H_{j,\nu}(x)+r_j(x),$ $\nu\in\mathbf{N}.$
Here $\{H_{j,\nu}\}$ is a sequence of polynomials converging
locally uniformly to $H_{j}$ on $\tilde{U}.$ It is clear that
$R(x,f^{\nu}(x))=W(x)g_{\nu}(x),$ for some
$g_{\nu}\in\mathcal{O}(\tilde{U}).$ The function $g$ is
non-vanishing on $\overline{U_0}$ therefore shrinking $\tilde{U},$
if necessary, we complete the proof for $n=1.$

Now suppose that $n>1.$ Let $f,R$ be a holomorphic mapping and a
polynomial respectively satisfying the hypotheses of the
proposition. Fix $U_0\in\mathcal{U}_{f,R}.$ By the definition of
$\mathcal{U}_{f,R}$ there exists a biholomorphism
$\phi:\check{U}\rightarrow\hat{U},$ where
$\check{U}\subset\mathbf{C}^n$ is a Runge domain,
$U_0\subset\subset\hat{U}$ and
$\check{U}\subset\check{U}_1\times\mathbf{C}$ for some open
connected $\check{U}_1\subset\mathbf{C}^{n-1},$ such that (1) and
(2) of Subsection \ref{construction} are satisfied.

Next observe that to complete the proof it is sufficient to show
that there is an open $E$ with $\phi^{-1}(U_0)\subset\subset
E\subset\check{U}$ and there are sequences $\{\tilde{f}^{\nu}\},$
$\{\phi^{\nu}\}$ of Nash mappings converging locally uniformly on
$E$ to $\tilde{f}=f\circ\phi,\phi$ respectively in such a way that
$\{\{x\in E :R(\phi^{\nu}(x),\tilde{f}^{\nu}(x))=0\}\}$ converges
in the sense of chains to $\{x\in E: R(\phi(x),\tilde{f}(x))=0\}.$
Indeed, given such sequences we may assume, shrinking $E$ if
necessary, that $\phi^{\nu}|_E$ is invertible for almost all
$\nu.$ Moreover, there is an open $\tilde{U}\subset\phi(E)$ such
that $U_0\subset\subset\tilde{U}\subset\phi^{\nu}(E)$ for almost
all $\nu.$ Consequently, $\{\{x\in\tilde{U}
:R(x,\tilde{f}^{\nu}\circ(\phi^{\nu})^{-1}(x))=0\}\}$ converges to
$\{x\in\tilde{U}: R(x,{f}(x))=0\}$ and we may set
$f^{\nu}=\tilde{f}^{\nu}\circ(\phi^{\nu})^{-1}.$

Before approximating $\tilde{f},\phi$ we show that there are Nash
mappings
$$A_l^{\nu}=(a^{\nu}_{l,1},\ldots,a^{\nu}_{l,p_l}),\mbox{ }
B_j^{\nu}=(b^{\nu}_{j,0},\ldots,b^{\nu}_{j,d-1}),\mbox{ }
C_i^{\nu}=(c^{\nu}_{i,0},\ldots,c^{\nu}_{i,d-1}),$$ for
$l=1,\ldots,\hat{s},$ $j=1,\ldots,k,$ $i=1,\ldots,n,$
$\nu\in\mathbf{N},$ defined on some open set
$E_1\subset\mathbf{C}^{n-1}$ with
$\pi(\phi^{-1}(U_0))\subset\subset E_1\subset\check{U}_1$ such
that the following hold. The sequence
$\{g^{\nu}:E_1\rightarrow\mathbf{C}^s\},$ where
$g^{\nu}=(A_1^{\nu},\ldots,A_{\hat{s}}^{\nu},B_1^{\nu},\ldots,B_{k}^{\nu},C_{1}^{\nu},\ldots,C_n^{\nu}),$
converges uniformly to $g|_{E_1}$ and $T_1\circ
g^{\nu}=\ldots=T_{d}\circ g^{\nu}=0$ for $\nu\in\mathbf{N}.$ Here
$g$ is the mapping from the condition (2) and $T_1,\ldots,T_{d}$
are polynomials given by the equation (*) of Subsection
\ref{construction}.

To this end, observe that by (2) it holds
$\pi(\phi^{-1}(U_0))\in\mathcal{U}_{g,S},$ where the polynomial
$S$ is described in the previous subsection. By the properties of
$S$ it is sufficient to show that there is a sequence
$\{h^{\nu}:\tilde{E}_1\rightarrow\mathbf{C}^s\}$ of Nash mappings
converging locally uniformly to $g|_{\tilde{E}_1},$ where
$\pi(\phi^{-1}(U_0))\subset\subset\tilde{E}_1\subset\check{U}_1,$
such that $\{\{x'\in\tilde{E}_1: S(x',h^{\nu}(x'))=0\}\}$
converges to $\{x'\in\tilde{E}_1: S(x',g(x'))=0\}$ in the sense of
chains (then $E_1$ may be taken to be any open set with
$\pi(\phi^{-1}(U_0))\subset\subset E_1\subset\subset\tilde{E}_1$).
This in turn is immediate by the induction hypothesis.

Using the components of $A_l^{\nu}, B_j^{\nu}, C_i^{\nu}$ define on $E=(E_1\times\mathbf{C})\cap\check{U}$ the following functions:\vspace*{2mm}\\
$W^{\nu}_l(x)=x_n^{p_l}+x_n^{p_l-1}a^{\nu}_{l,1}(x')+\ldots+a^{\nu}_{l,p_l}(x'),$\\
$r^{\nu}_j(x)=x_n^{d-1}b^{\nu}_{j,0}(x')+x_n^{d-2}b^{\nu}_{j,1}(x')+\ldots+b^{\nu}_{j,d-1}(x'),$\\
$\check{r}^{\nu}_i(x)=x_n^{d-1}c^{\nu}_{i,0}(x')+x_n^{d-2}c^{\nu}_{i,1}(x')+\ldots+c^{\nu}_{i,d-1}(x'),$\vspace*{2mm}\\
for $l=1\ldots,\hat{s}, j=1,\ldots,k, i=1\ldots,n.$ Next put
$W^{\nu}=(W^{\nu}_1)^{k_1}\cdot\ldots\cdot
(W^{\nu}_{\hat{s}})^{k_{\hat{s}}},$ where $k_j$ is the
multiplicity of the factor $W_j$ of $W$ (see
Subsection~\ref{construction}). Now define
$\tilde{f}^{\nu}=(f^{\nu}_1,\ldots,f^{\nu}_k),$
$\phi^{\nu}=(\phi^{\nu}_1,\ldots,\phi^{\nu}_n)$ by
$$f^{\nu}_j=W^{\nu}H^{\nu}_{j}+r^{\nu}_j,\mbox{ }
\phi^{\nu}_i=W^{\nu}\check{H}^{\nu}_{i}+\check{r}^{\nu}_i$$ for
$j=1,\ldots,k,$ $i=1,\ldots,n.$ Here
$\{H^{\nu}_{j}\},\{\check{H}^{\nu}_{i}\},$ are any sequences of
polynomials converging locally uniformly on $E$ to
$H_{j},\check{H}_{i}$ respectively. (Recall that
$H_{j},\check{H}_{i}$ are obtained in Subsection
\ref{construction} dividing $f_j, \phi_i$ by $W.$ The existence of
$\{H^{\nu}_{j}\},\{\check{H}^{\nu}_{i}\}$ follows by the fact that
$\check{U}$ is a Runge domain.) Clearly, $\{\tilde{f}^{\nu}\},
\{\phi^{\nu}\}$ converge locally uniformly to
$\tilde{f}|_{E},\phi|_{E}$ respectively.

Finally, the equation (*) from Subsection \ref{construction}, in
view of the fact that $T_1\circ g^{\nu}=\ldots=T_{d}\circ
g^{\nu}=0,$ implies
$R(\phi^{\nu}(x),\tilde{f}^{\nu}(x))=\tilde{H}^{\nu}(x)W^{\nu}(x)$
for every $x\in E,$ where $\tilde{H}^{\nu}\in\mathcal{O}(E).$
Since $\{W_l^{\nu}|_E\}$ converges to ${W_l}|_E$ locally
uniformly, for $l=1,\ldots,\hat{s}$ (where
$W_1,\ldots,W_{\hat{s}}$ are optimal polynomials such that
$W=W_1^{k_1}\cdot\ldots\cdot W_{\hat{s}}^{k_{\hat{s}}}$ and
$dim(W^{-1}_i(0)\cap W^{-1}_j(0))<n-1$ for every $i\neq j$) it
holds: $\{\{x\in E: W^{\nu}(x)=0\}\}$ converges to $\{x\in E:
W(x)=0\}$ in the sense of chains. The function $\tilde{H}$ given
by (1) is non-vanishing on $\overline{\phi^{-1}(U_0)}$ therefore
shrinking $E$ if necessary we obtain the required claim.\qed
%
%
\subsubsection{Algorithm}\label{examp} Based
on the proof of Theorem \ref{semiglobmain}, we present a recursive
algorithm of Nash appro\-ximation of a holomorphic mapping
$f:U\rightarrow V\subset\mathbf{C}^{\hat{m}},$ where $U$ is a
domain in $\mathbf{C}^n$ and $V$ is an algebraic variety. For
$\nu\in\mathbf{N},$ the approximating mapping
$f^{\nu}=(f^{\nu}_1,\ldots,f^{\nu}_{\hat{m}}):U_0\rightarrow V,$
returned as the output of the algorithm, is represented by
$\hat{m}$ non-zero polyno\-mials
$P_i^{\nu}(x,z_i)\in(\mathbf{C}[x])[z_i],$ $i=1,\ldots,\hat{m},$
such that $P_i^{\nu}(x,f^{\nu}_i(x))=0$ for $x\in U_0.$ For
simplicity we restrict attention to the local case i.e. $U_0$ is
an open neighborhood of a fixed $x_0\in U.$
More precisely, we work with the following data:\vspace*{2mm}\\
\textbf{Input:} a holomorphic mapping
$f=(f_1,\ldots,f_{\hat{m}}):U\rightarrow
V\subset\mathbf{C}^{\hat{m}},$ $f=f(x),$ where $U$ is an open
neighborhood of $0\in\mathbf{C}^n$ and $V$ is an algebraic
variety.\vspace*{2mm}\\
\textbf{Output:} $P_i^{\nu}(x,z_i)\in(\mathbf{C}[x])[z_i],$
$P_i^{\nu}\neq 0$ for $i=1,\ldots,\hat{m}$ and $\nu\in\mathbf{N},$
with the following properties:\vspace*{1mm}\\
(a) $P_i^{\nu}(x,f_i^{\nu}(x))=0$ for every $x\in U_0,$ where
$f^{\nu}=(f^{\nu}_1,\ldots,f^{\nu}_{\hat{m}}):U_0\rightarrow V$ is
a holomorphic mapping such that $\{f^{\nu}\}$ converges uniformly
to $f$ on an open neighborhood $U_0$ of
$0\in\mathbf{C}^n,$\\
(b) $P_i^{\nu}$ is a unitary polynomial in $z_i$ of degree
independent of $\nu$ whose co\-efficients (belonging to
$\mathbf{C}[x]$) converge uniformly to holomorphic functions on
$U_0$ as $\nu$ tends to infinity.\vspace*{2mm}

Before going into detail let us comment on the notation and the
idea of the algorithm. First, the meaning of the symbols
$V_{(m)},$ $\mathcal{S}(V)$ and the notion of the optimal
polynomial used below can be found in Subsection
\ref{setswithprop}. Next, in steps 2 and 5 we apply linear changes
of the coordinates. Having approximated the mapping $\hat{J}\circ
f\circ J|_{J^{-1}(U)}:J^{-1}(U)\rightarrow\hat{J}(V),$ where
$\hat{J}:\mathbf{C}^{\hat{m}}\rightarrow\mathbf{C}^{\hat{m}},$
$J:\mathbf{C}^n\rightarrow\mathbf{C}^n$ are linear isomorphisms,
one can obtain the output data for $f$ following standard
arguments. (Composing $f$ and $J$ does not lead to any
difficulties. As for $\hat{J},$ it is sufficient to use the fact
that the integral closure of a commutative ring in another
commutative ring is again a ring.) Therefore, when the coordinates
are changed, we write what (as a result) may be assumed about the
mapping $f,$ but the notation is left unchanged.

The aim of steps 1-3 is to prepare the variety $V$ so that the
polynomial $R$ calculated in step 4 satisfies the assertion of
Proposition \ref{claglo} (cp. the proof of Proposition
\ref{claglo}). Steps 5-9 are responsible for the fact that for
$f^{\nu}_1,\ldots,f^{\nu}_m$ defined in step 10 the sequence
$\{\{R(f^{\nu}_1(x),\ldots,f^{\nu}_m(x))=0\}\}$ converges to
$\{R(f_1(x),\ldots,f_m(x))=0\}$ in the sense of chains, in a
neighborhood of $0\in\mathbf{C}^n,$ as $\nu$ tends to infinity.
This property implies (cp. the proof Proposition \ref{claglo})
that there is an open neighborhood $U_0$ of $0\in\mathbf{C}^n$
such that for $\nu$ large enough the set
$\{(x,z_{m+1},\ldots,z_{m+s})\in U_0\times\mathbf{C}^s:
(f^{\nu}_1(x),\ldots,f^{\nu}_m(x),z_{m+1},\ldots,z_{m+s})\in V\}$
contains a graph of the mapping
$x\mapsto(f^{\nu}_{m+1}(x),\ldots,f^{\nu}_{m+s}(x))$ approximating
the mapping $x\mapsto(f_{m+1}(x),\ldots,f_{m+s}(x))$ (here
$\hat{m}=m+s$). The latter fact is used in step 11 to calculate
$P^{\nu}_{m+1},\ldots,P^{\nu}_{m+s}.$ As for
$P^{\nu}_{1},\ldots,P^{\nu}_{m},$ these polynomials are obtained
in step 10 by applying the results of the algorithm switched for
the lower dimensional case in step 9.\vspace*{2mm}\\
\textbf{Algorithm:} \textbf{1.} If $f(U)\subseteq Sing(V)$ then
repeat replacing $V$ by $Sing(V)$ until $f(U)\nsubseteq Sing(V).$
Next replace $V$ by $V_{(m)}$ such that
$f(U)\subset V_{(m)}.$\\
\textbf{2.} Apply a linear change of the coordinates in
$\mathbf{C}^{\hat{m}}$ after which $\rho|_{V}$ is a proper mapping
and $\rho(f(U))\nsubseteq\rho(\mathcal{S}(V)),$ where
$\rho:\mathbf{C}^m\times\mathbf{C}^s\approx\mathbf{C}^{\hat{m}}\rightarrow\mathbf{C}^m$
is the natural projection.\\
\textbf{3.} Choose a $\mathbf{C}$-linear form
$L:\mathbf{C}^s\rightarrow\mathbf{C}$ such that the generic fibers
of $\rho|_{V}$ and $\tilde{\rho}|_{\Phi_L(V)}$ over $\mathbf{C}^m$
have the same cardinalities and
$\rho(f(U))\nsubseteq\tilde{\rho}(\mathcal{S}(\Phi_L(V))).$ Here
$\tilde{\rho}:\mathbf{C}^m\times\mathbf{C}\rightarrow\mathbf{C}^m$
is the natural projection and $\Phi_L(y,v)=(y,L(v))$ for
$(y,v)\in\mathbf{C}^m\times\mathbf{C}^s.$\\
\textbf{4.} Calculate the discriminant $R\in\mathbf{C}[y]$ of the
optimal polynomial $P(y,z)\in(\mathbf{C}[y])[z]$ describing
$\Phi_L(V)\subset\mathbf{C}^m_y\times\mathbf{C}_z.$\\
\textbf{5.} Apply a linear change of the coordinates in
$\mathbf{C}^n$ after which $R(\rho(f(x)))=\tilde{H}(x)W(x)$ in
some neighborhood of $0\in\mathbf{C}^n,$ where $\tilde{H}$ is a
holomorphic function, $\tilde{H}(0)\neq 0$ and $W$ is a unitary
polynomial in $x_n$ with holomorphic coefficients\linebreak
depending on $x'=(x_1,\ldots,x_{n-1})$ each of which vanishes at
$0\in\mathbf{C}^{n-1}.$ Put $d=deg(W).$\\
\textbf{6.} Divide $f_i$ by $W$ to obtain
$f_i(x)=W(x)H_i(x)+r_i(x)$ in some neighborhood of
$0\in\mathbf{C}^n,$ $i=1,\ldots,m.$ Here $H_i$ is a holomorphic
function and $r_i$ is a polynomial in
$x_n,$ $deg(r_i)<d,$ with holomorphic coefficients depending on $x'.$\\
\textbf{7.} Find optimal polynomials $W_1,\ldots,W_{\hat{s}}$ in
$x_n$ with holomorphic coefficients depending on $x'$ such that
$W=W_1^{k_1}\cdot\ldots\cdot W_{\hat{s}}^{k_{\hat{s}}}$ and
$dim(W^{-1}_i(0)\cap W^{-1}_j(0))<n-1$ for every $i\neq j.$\\
\textbf{8.} Treating $H_i,$ $i=1,\ldots,m,$ and all the
coefficients of $W_1,\ldots, W_{\hat{s}},r_1,\ldots,r_m$ as new
variables (except for the coefficient $1$ standing at the leading
terms of $W_1,\ldots,W_{\hat{s}}$) apply the division procedure
for
polynomials to obtain:\\
$R(WH_1+r_1,\ldots,WH_m+r_m)=\tilde{W}W+x_n^{d-1}T_1+x_n^{d-2}T_2+\ldots+T_{d}.$
Here $T_1,\ldots,T_{d}$ are polynomials depending only on the
variables standing for the coefficients of
$W_1,\ldots,W_{\hat{s}},$ $r_1,\ldots,r_m.$ Moreover,
$T_1(g)=\ldots=T_{d}(g)=0,$ where $g$ is the holomorphic mapping
whose
components are these coefficients (cp. Subsection \ref{construction}).\\
\textbf{9.} If $g$ is not constant (i.e. it depends on $n-1\geq 1$
variables) then apply\linebreak the Algorithm with $f, V$ replaced
by $g$ and $ \{u\in\mathbf{C}^{\hat{d}}:
T_1(u)=\ldots=T_{d}(u)=0\}$ respectively, where $\hat{d}$ is the
number of the components of $g.$ As a result, for\linebreak every
$c(x')$ which is a coefficient of some of
$W_1,\ldots,W_{\hat{s}},$ $r_1,\ldots,r_m$ one obtains a sequence
$\{Q_c^{\nu}(x',t_c)\}$ of unitary polyno\-mials satisfying (a)
and (b) above with $x,$ $z_i,$ $\{f^{\nu}\}$ replaced by $x',$
$t_c,$ $\{g^{\nu}\}$ respectively. Here $\{g^{\nu}\}$ is a
sequence of Nash mappings converging to $g$ in some neighborhood
of $0\in\mathbf{C}^{n-1}$ such that $T_1\circ
g^{\nu}=\ldots=T_d\circ g^{\nu}=0$ for every $\nu\in\mathbf{N}$
(cp. Subsection \ref{globperties}). If $g$ is
constant then it is its own approximation yielding the $Q^{\nu}_c$'s immediately.\\
\textbf{10.} Approximate $H_i,$ for $i=1,\ldots,m,$ by a sequence
$\{H_i^{\nu}\}$ of polynomials.\linebreak Let
$W_1^{\nu},\ldots,W_{\hat{s}}^{\nu},r_1^{\nu},$
$\ldots,r_m^{\nu},$ for every $\nu\in\mathbf{N},$ be the
polynomials in $x_n$ defined by replacing the coefficients of
$W_1,\ldots,W_{\hat{s}},$ $r_1,\ldots,r_m$ by their Nash
approxima\-tions (i.e. the components of $g^{\nu}$) determined in
step 9. Using $Q_c^{\nu}$ (for all $c$) and $H^{\nu}_i$ one can
calculate $P_i^{\nu}\in(\mathbf{C}[x])[z_i],$ for $i=1,\ldots,m,$
satisfying (b) and (a) with
$f_i^{\nu}=H_i^{\nu}(W_1^{\nu})^{k_1}\cdot\ldots\cdot(W_{\hat{s}}^{\nu})^{k_{\hat{s}}}+r_i^{\nu}$
being the $i$'th component of the mapping $f^{\nu}$ (whose last
$\hat{m}-m$ components are determined by
$P_{m+1}^{\nu},\ldots,P_{\hat{m}}^{\nu}$ obtained in the next
step). To calculate $P^{\nu}_1,\ldots,P^{\nu}_m$ one can follow
the standard proof of the fact that the integral closure of a
commutative ring in another commutative ring is again a ring.\\
\textbf{11.} Put
$V^{\nu}=\{(x,z)\in\mathbf{C}^n_x\times\mathbf{C}^{m+s}_{z} : z\in
V, P^{\nu}_i(x,z_i)=0 \mbox{ for } i=1,\ldots,m\},$ where
$z=(z_1,\ldots,z_m,z_{m+1},\ldots,z_{m+s}).$ For $i=1,\ldots,s$
and $\nu\in\mathbf{N}$ take
$P^{\nu}_{m+i}\in(\mathbf{C}[x])[z_{m+i}]$ to be the optimal
polynomial describing the image of the projection of $V^{\nu}$
onto $\mathbf{C}^n_x\times\mathbf{C}_{z_{m+i}}.$


\end{document}